\documentclass[11pt]{article}
\usepackage{amsmath}
\usepackage{amsfonts}
\usepackage{graphicx}
\usepackage{amssymb}
\usepackage{cite}

\newtheorem{condition**}{A*}
\newtheorem{condition***}{C*}
\newtheorem{condition*}{C}

\newtheorem{proposition}{Proposition}[section]

\newtheorem{definition}{Definition}[section]
\newtheorem{theorem}{Theorem}[section]
\newtheorem{lemma}{Lemma}[section]
\newtheorem{remark}{Remark}[section]

\newenvironment{keywords}{{\bf Key words: }}{}

\begin{document}

\title{Mean Field Linear-Quadratic-Gaussian (LQG) Games: Major and Minor Players}

\author{Jianhui Huang\thanks{J. Huang is with the Department of Applied Mathematics, The Hong Kong Polytechnic University, Hong Kong
(majhuang@polyu.edu.hk).}
\quad Shujun Wang \thanks{S. Wang is with the Department of Applied Mathematics, The Hong Kong Polytechnic University, Hong Kong
(shujun.wang@connect.polyu.hk).}
\quad Zhen Wu \thanks{Z. Wu is with School of Mathematics, Shandong University, Jinan 250100, China (wuzhen@sdu.edu.cn).}}

\maketitle

\begin{abstract}
This paper is concerned with a backward-forward stochastic differential equation (BFSDE) system, in which a large number of negligible agents are coupled in their dynamics via state average. Here some BSDE is introduced as the dynamics of major player, while dynamics of minor players are described by
SDEs. Some auxiliary mean-field SDEs (MFSDEs) and a $3\times2$ mixed forward-backward stochastic differential equation (FBSDE) system are considered and analyzed instead of involving the fixed-point analysis as in other mean-field games. We also derive the decentralized strategies which are shown to satisfy the $\epsilon$-Nash equilibrium property.
\end{abstract}


%%%%%%%%%%%%%%%%%%%%%%%%%%%%%%%%%%%%%%%%%%%%%%%%%%%%%%%%%%%%%%%%%%%%%%%% Section 1 %%%%%%%%%%%%%%%%%%%%%%%%%%%%%%%%%%%%%%%%%%%%%%%%%%%%%%%%%%%%%%%
\section{Introduction}
In recent years, within the context of noncooperative game theory, the dynamic optimization or control of stochastic large-population (LP) system has attracted consistent and intense attentions in a variety of fields including biology, engineering, operational research, mathematical finance and economics, social science, etc. The most special feature of controlled LP system lies in the existence of considerable insignificant agents whose dynamics and (or) cost functionals are coupled via the state-average across the whole population. To design low-complexity strategies, one efficient methodology is the associated mean-field games which enable us to obtain the decentralized control. Readers may refer \cite{ll} for the motivation and methodology, and \cite{BCQ} for recent progress in mean-field game theory. Besides, some other recent literature include \cite{B12,B11,dh,hcm07,HCM12,hmc06,LZ08} for linear-quadratic-Gaussian (LQG) mean-field games of large-population system, \cite{TZB11} for risk sensitive mean-field games. In contrast to previous work, in \cite{H10} M. Huang discussed large population systems with major and minor players by analyzing the case in an infinite set where the minor players are from a total of $K$ classes. Later on, Nguyen and Huang \cite{NH12} considered an LQ problem by directly treating the mean field $z$ in the population limit as a random process with random coefficients. In addition, J. Yong \cite{Yong02} investigated a leader-follower hierarchical game with open-loop information.

In reality, major player (the government, for example) may expect to obtain some recursive utilities who also has a significant role in affecting all minor players. For this kind of model, some backward stochastic differential equation (BSDE) is introduced as the dynamics of major player. Bismut \cite{B78} introduced BSDEs as the adjoint equations. Pardoux and Peng first proved the existence and uniqueness of solution for nonlinear BSDEs in \cite{PP}, which has been extensively used in stochastic control and mathematical finance. Independently, Duffie and Epstein \cite{DE} introduced BSDEs under economic background. Then El Karoui, Peng and Quenez gave the formulation of recursive utilities from the BSDE point of view. As found by \cite{EPQ}, the recursive utility process can be regarded as a solution of BSDE.

Furthermore, with common forward dynamics of minor players, in this paper we derive a backward-forward stochastic differential equation (BFSDE) system, in which a large number of negligible agents are coupled in their dynamics via state average. We discuss the related mean-field LQG games and derive the decentralized strategies which are shown to satisfy the $\epsilon$-Nash equilibrium property. A stochastic process which relates to the state of major player is introduced here to be the approximation of the state-average process. Some auxiliary mean-field stochastic differential equations (MFSDEs) and a $3\times2$ mixed forward-backward stochastic differential equation (FBSDE) system are considered and analyzed.

The rest of this paper is organized as follows. Section 2 formulates the large population LQG games of backward-forward systems. In Section 3, we derive the optimal control of auxiliary track system and the consistency conditions. Section 4 is devoted to the related $\epsilon$-Nash equilibrium. Section 5 concludes our work.

%%%%%%%%%%%%%%%%%%%%%%%%%%%%%%%%%%%%%%%%%%%%%%%%%%%%%%%%%%%%%%%%%%%%%%%% Section 2 %%%%%%%%%%%%%%%%%%%%%%%%%%%%%%%%%%%%%%%%%%%%%%%%%%%%%%%%%%%%%%%
\section{Problem formulation}

Throughout this paper, we denote by $\mathbb{R}^m$ the $m$-dimensional Euclidean space. Consider a finite time horizon $[0,T]$ for a fixed $T>0$. Suppose $(\Omega,
\mathcal F, \{\mathcal F_t\}_{0\leq t\leq T}, P)$ is a complete
filtered probability space on which a standard $(d+m\times N)$-dimensional Brownian motion $\{W_0(t),W_i(t),\ 1\le i\leq N\}_{0 \leq t \leq T}$ is defined.
$\mathcal F^{w_0}_t:=\sigma\{W_0(s), 0\leq s\leq t\}$, $\mathcal F^{w_i}_t:=\sigma\{W_i(s), 0\leq s\leq t\}$, $\mathcal F^{i}_t:=\sigma\{W_0(s),W_i(s);0\leq s\leq t\}$. Here, $\{\mathcal F^{w_0}_t\}_{0\leq t\leq T}$ stands for the information of the major player; while $\{\mathcal F^{w_i}_t\}_{0\leq t\leq T}$ the individual information of $i^{th}$ minor player. For a given filtration $\{\mathcal G_t\}_{0\leq t\leq T},$ let $L^{2}_{\mathcal{G}_t}(0, T; \mathbb{R}^m)$ denote the space of all $\mathcal{G}_t$-progressively measurable processes with values in $\mathbb{R}^m$ satisfying $\mathbb{E}\int_0^{T}|x(t)|^{2}dt<+\infty;$ $L^{2}(0, T; \mathbb{R}^m)$ the space of all deterministic functions defined on $[0,T]$ in $\mathbb{R}^m$ satisfying $\int_0^{T}|x(t)|^{2}dt<+\infty;$ $C(0,T;\mathbb{R}^m)$ the space of all continuous functions defined on $[0,T]$ in $\mathbb{R}^m$. For simplicity, in what follows we focus on the $1$-dimensional processes, which means $d=m=1$.

Consider an large population system with $(1+N)$ individual agents, denoted by $\mathcal{A}_{0}$ and $\{\mathcal{A}_{i}\}_{1 \leq i \leq N},$ where $\mathcal{A}_{0}$ stands for the major player, while $\mathcal{A}_{i}$ stands for $i^{th}$ minor player. The dynamics of $\mathcal{A}_{0}$ is given by a BSDE as follows:
\begin{equation}\label{e1}
    \left\{
    \begin{aligned}
    dx_0(t)=& \big[A_0x_{0}(t)+B_0u_{0}(t)+C_0z_0(t)\big]dt+z_0(t)dW_{0}(t),\\
    x_{0}(T)=& \xi
    \end{aligned}
    \right.
\end{equation}
where $\xi\in \mathcal{F}^{w_0}_T$ satisfies $E|\xi|^2<+\infty.$ The state of minor player $\mathcal{A}_{i}$ is a SDE satisfying
\begin{equation}\label{e2}
    \left\{
    \begin{aligned}
    dx_{i}(t)=& \big[Ax_{i}(t)+Bu_{i}(t)+Dx^{(N)}(t)+\alpha x_0(t)\big]dt+\sigma dW_{i}(t),\\
    x_{i}(0)=& x_{i0}
    \end{aligned}
    \right.
\end{equation}
where $x^{(N)}(t)=\frac{1}{N}\sum\limits_{i=1}^Nx_i(t)$ is the state-average of minor players; $x_{i0}$ is the initial value of $\mathcal{A}_i$. Here, $A_0,B_0,C_0,A, B, D, \alpha, \sigma$ are scalar constants. Assume that $\mathcal F_t$ is the augmentation of $\sigma\{W_0(s),W_i(s),x_{i0};0\leq s\leq t,1\leq i\leq N\}$ by all the $P$-null sets of $\mathcal{F}$, which is the full information accessible to the large population system up to time $t$. The admissible control strategy $u_0\in \mathcal{U}_0,u_i\in \mathcal{U}_i$ where$$
\mathcal{U}_0\triangleq\big\{u_0|u_0(t)\in L^{2}_{\mathcal{F}^{w_0}_t}(0, T; \mathbb{R})\big\},
$$
and $$
\mathcal{U}_i\triangleq\big\{u_i|u_i(t)\in L^{2}_{\mathcal{F}_t}(0, T; \mathbb{R})\big\},\ 1\leq i \leq N.
$$
Let $u=(u_0,u_1, \cdots, u_{N})$ denote the set of control strategies of all $(1+N)$ agents; $u_{-i}=(u_0,u_1, \cdots, u_{i-1},$ $u_{i+1}, \cdots u_{N})$ the control strategies except $i^{th}$ agent $\mathcal{A}_i,0\leq i\leq N.$ The cost functional for $\mathcal{A}_0$ is given by
\begin{equation}\label{e3}
    \begin{aligned}
    J_{0}(u_{0}(\cdot), u_{-0}(\cdot))=&\frac{1}{2} \mathbb{E}\left\{ \int_{0}^{T}\left[Q_0\Big(x_{0}(t)-x^{(N)}(t)\Big)^{2}+R_0u_{0}^{2}(t)\right]dt+H_{0}x_{0}^{2}(0) \right\}
\end{aligned}
\end{equation}
where $Q_0\geq 0, R_0>0, H_{0}\geq 0$. The individual cost functional for $\mathcal{A}_i,1\leq i\leq N$ is
\begin{equation}\label{e4}
    \begin{aligned}
    J_{i}(u_{i}(\cdot), u_{-i}(\cdot))=&\frac{1}{2} \mathbb{E}\left\{ \int_{0}^{T}\left[Q\Big(x_{i}(t)-x^{(N)}(t)\Big)^{2}+Ru_{i}^{2}(t)\right]dt+Hx_{i}^{2}(T) \right\}
\end{aligned}
\end{equation}
where $Q\geq 0, R>0, H\geq 0.$

\begin{remark}\emph{Unlike the systems introduced in \cite{H10, NH12}, the dynamics of major player is a BSDE. The term $H_0x_{0}^{2}(0)$ is introduced in \eqref{e3} to denote some recursive evaluation or nonlinear expectation. One practical meaning is the initial hedging deposits in the pension fund industry. For simplicity, behaviors of major player (the government, as an example) affect the state of each minor player (the peasant); while the state-average of minor players appears only in the cost functional of major player.} \end{remark}

We introduce the following assumption:
\begin{description}
  \item[(H1)] $\{x_{i0}\}_{i=1}^N$ are independent and identically distributed (i.i.d) with $\mathbb{E}x_{i0}=x$, $\mathbb{E}|x_{i0}|^2<+\infty,$ and also independent of $\{W_0,W_i,1\leq i\leq N\}$.
\end{description}
 Now, we formulate the large population dynamic optimization problem.\\

\textbf{Problem (I).}
Find a control strategies set $\bar{u}=(\bar{u}_0,\bar{u}_1,\cdots,\bar{u}_N)$ which satisfies
$$
J_i(\bar{u}_i(\cdot),\bar{u}_{-i}(\cdot))=\inf_{u_i\in \mathcal{U}_i}J_i(u_i(\cdot),\bar{u}_{-i}(\cdot))
$$where $\bar{u}_{-i}$ represents $(\bar{u}_0,\bar{u}_1,\cdots,\bar{u}_{i-1},\bar{u}_{i+1},\cdots, \bar{u}_N),0\leq i\leq N$.

%%%%%%%%%%%%%%%%%%%%%%%%%%%%%%%%%%%%%%%%%%%%%%%%%%%%%%%%%%%%%%%%%%%%%%%% Section 3 %%%%%%%%%%%%%%%%%%%%%%%%%%%%%%%%%%%%%%%%%%%%%%%%%%%%%%%%%%%%%%%
\section{The Limiting Optimal Control and NCE Equation System}
To study Problem \textbf{(I)}, one efficient approach is to discuss the associated mean-field games via limiting problem when the agent number $N$ tends to infinity. For a given stochastic triple $(\bar{x}(\cdot),k(\cdot),\beta_0(\cdot))\in L^{2}_{\mathcal{F}^{w_0}_t}(0, T; \mathbb{R})$, suppose $x^{(N)}(\cdot)$ is approximated by $\bar{x}(\cdot)$ as $N \rightarrow +\infty.$

Introduce the following auxiliary dynamics of major and minor players, still denote by $x_0(\cdot),x_i(\cdot)$ respectively:
\begin{equation}\label{e5}
    \left\{
    \begin{aligned}
    &dx_0(t)= \big[A_0x_{0}(t)+B_0u_{0}(t)+C_0z_0(t)\big]dt+z_0(t)dW_{0}(t),\\
    &x_{0}(T)= \xi,\\
    &d\bar{x}(t)=\big[\bar{A}\bar{x}(t)+\bar{B}x_0(t)+\bar{C}k(t)\big]dt,\\
    &\bar{x}(0)=x,\\
    &dk(t)=\big[\tilde{A}k(t)+\tilde{B}\bar{x}(t)+\tilde{C}x_0(t)\big]dt+\beta_0(t)dW_0(t),\\
    &k(T)=0
    \end{aligned}
    \right.
\end{equation}
and
\begin{equation}\label{e6}
    \left\{
    \begin{aligned}
    dx_{i}(t)=& \big[Ax_{i}(t)+Bu_{i}(t)+D\bar{x}(t)+\alpha x_0(t)\big]dt+\sigma dW_{i}(t),\\
    x_{i}(0)=& x_{i0}.
    \end{aligned}
    \right.
\end{equation}
Here $\bar{\varphi}$ and $\tilde{\varphi},\varphi=A,B,C$ are to be determined. The associated limiting cost functionals become
\begin{equation}\label{e7}
    \begin{aligned}
    \bar{J}_{0}(u_{0}(\cdot))=&\frac{1}{2} \mathbb{E}\left\{ \int_{0}^{T}\left[Q_0\Big(x_{0}(t)-\bar{x}(t)\Big)^{2}+R_0u_{0}^{2}(t)\right]dt+H_{0}x_{0}^{2}(0) \right\}
\end{aligned}
\end{equation}
and
\begin{equation}\label{e8}
    \begin{aligned}
    \bar{J}_{i}(u_{i}(\cdot))=&\frac{1}{2} \mathbb{E}\left\{ \int_{0}^{T}\left[Q\Big(x_{i}(t)-\bar{x}(t)\Big)^{2}+Ru_{i}^{2}(t)\right]dt+Hx_{i}^{2}(T) \right\}.
\end{aligned}
\end{equation}
\begin{remark}\emph{Due to the reason that the state-average of minor players appears only in the cost functional of major player, the first equation in \eqref{e5} has the same form with \eqref{e1} actually. However, for regularity, we still write it out.   } \end{remark}

Thus, we formulate the limiting LQG game \textbf{(II)} as follows. \\

\textbf{Problem (II).} For $i^{th}$ agent $\mathcal{A}_i$, $i=0,1,2,\cdots,N,$ find $\bar{u}_i\in \mathcal{U}_i$ satisfying
\begin{equation}\label{e9}
\bar{J}_i(\bar{u}_i(\cdot))=\inf_{u_i\in \mathcal{U}_i}\bar{J}_i(u_i(\cdot)).
\end{equation}
$\bar{u}_i$ satisfying \eqref{e9} is called an optimal control for (\textbf{II}).

To get the optimal control of Problem (\textbf{II}), we should obtain the optimal control of $\mathcal{A}_0$ first. For this aim, we have the following lemma.

\begin{lemma}\label{l1}
Corresponding to the forward-backward system \eqref{e5} and \eqref{e7}, the optimal control of $\mathcal{A}_0$ for \emph{(\textbf{II})} is given by
\begin{equation}\label{e10}
    \bar{u}_{0}(t)=-B_0R_0^{-1}p_0(t)
\end{equation}where the adjoint process $p_0(\cdot)$ and the corresponding optimal trajectory $(\hat{x}_0(\cdot),\hat{z}_0(\cdot))$ satisfy the following Hamilton system
 \begin{equation}\label{e11}
    \left\{
    \begin{aligned}
    &d\hat{x}_0(t)= \big[A_0\hat{x}_0(t)-B_0^2R_0^{-1}p_0(t)+C_0\hat{z}_0(t)\big]dt+\hat{z}_0(t)dW_{0}(t),\\
    &d\bar{x}(t)=\big[\bar{A}\bar{x}(t)+\bar{B}\hat{x}_0(t)+\bar{C}k(t)\big]dt,\\
    &dk(t)=\big[\tilde{A}k(t)+\tilde{B}\bar{x}(t)+\tilde{C}\hat{x}_0(t)\big]dt+\beta(t)dW_0(t),\\
    &dp_0(t)=\big[-A_0p_0(t)-Q_0(\hat{x}_0(t)-\bar{x}(t))-\bar{B}p(t)-\tilde{C}q(t)\big]dt-C_0p_0(t)dW_0(t),\\
    &dp(t)=\big[-\bar{A}p(t)+Q_0(\hat{x}_0(t)-\bar{x}(t))-\tilde{B}q(t)\big]dt+\bar{\beta}(t)dW_0(t),\\
    &dq(t)=\big(-\tilde{A}q(t)-\bar{C}p(t)\big)dt,\\
    &\hat{x}_{0}(T)= \xi,\  \bar{x}(0)=x,\  k(T)=0,\   p_0(0)=-H_0\hat{x}_0(0),\  p(T)=0, \ q(0)=0
    \end{aligned}
    \right.
\end{equation}where $\beta(\cdot),\bar{\beta}(\cdot)\in L^{2}_{\mathcal{F}^{w_0}_t}(0, T; \mathbb{R})$.
\end{lemma}
{\it Proof.} For the variation of control $\delta u_0(\cdot)\in L^{2}_{\mathcal{F}^{w_0}_t}(0, T; \mathbb{R})$, introduce the following variational equations:
 \begin{equation}\label{e12}
    \left\{
    \begin{aligned}
    &d\delta x_0(t)= \big[A_0\delta x_0(t)+B_0 \delta u_0(t)+C_0\delta z_0(t)\big]dt+\delta z_0(t)dW_{0}(t),\\
    &d\delta\bar{x}(t)=\big[\bar{A}\delta\bar{x}(t)+\bar{B}\delta x_0(t)+\bar{C}\delta k(t)\big]dt,\\
    &d\delta k(t)=\big[\tilde{A}\delta k(t)+\tilde{B}\delta\bar{x}(t)+\tilde{C}\delta x_0(t)\big]dt+\delta\beta_0(t)dW_0(t),\\
    &\delta x_{0}(T)= 0,\  \delta\bar{x}(0)=0,\  \delta k(T)=0.
    \end{aligned}
    \right.
\end{equation}
Applying It\^o's formula to $p_0(t)\delta x_0(t)+p(t) \delta \bar{x}(t)+q(t)\delta k(t)$ and noting the associated first order variation of cost functional :
$$0=\delta\bar{J}_{0}(\bar{u}_{0})=\mathbb{E}\left\{ \int_{0}^{T}\left[Q_0\big(\hat{x}_{0}(t)-\bar{x}(t)\big)\big(\delta x_{0}(t)-\delta\bar{x}(t)\big)+R_0\bar{u}_{0}(t)\delta u_0(t)\right]dt+H_{0}\hat{x}_{0}(0)\delta x_0(0) \right\},$$
we obtain the optimal control \eqref{e10}. Combining all state and adjoint equations, and applying $\bar{u}_0(\cdot)$ to $\mathcal{A}_0$, we get the Hamilton system \eqref{e11}.
   \hfill   $\Box$

After obtaining the optimal control of major player $\mathcal{A}_0$, we aim to get the optimal control of minor player $\mathcal{A}_i$ corresponding to $\hat{x}_0(\cdot)$ respectively.
\begin{lemma}\label{l2}
Under \emph{(H1)}, the optimal control of $\mathcal{A}_i$ for \emph{(\textbf{II})} is
\begin{equation}\label{e13}
    \bar{u}_{i}(t)=-BR^{-1}p_i(t)
\end{equation}where the adjoint process $p_i(\cdot)$ and the corresponding optimal trajectory $\hat{x}_i(\cdot)$ satisfy the BSDE
 \begin{equation}\label{e14}
    \left\{
    \begin{aligned}
    dp_i(t)=& \big[-Ap_i(t)-Q\big(\hat{x}_i(t)-\bar{x}(t)\big)\big]dt+\beta_0(t)dW_0(t)+\beta_i(t)dW_i(t),\\
        p_i(T)=& H\hat{x}_i(T)
    \end{aligned}
    \right.
\end{equation}and SDE
\begin{equation}\label{e15}
    \left\{
    \begin{aligned}
    d\hat{x}_i(t)=& \big[A\hat{x}_i(t)-B^2R^{-1}p_i(t)+D\bar{x}(t)+\alpha \hat{x}_0(t)\big]dt+\sigma(t)dW_{i}(t),\\
    \hat{x}_{i}(0)=& x_{i0}.
    \end{aligned}
    \right.
\end{equation}
\end{lemma}
Here $\beta_0(\cdot),\beta_i(\cdot)\in L^{2}_{\mathcal{F}^i_t}(0, T; \mathbb{R})$; $\hat{x}_0(\cdot)$ and $\bar{x}(\cdot)$ are given by \eqref{e11}. For the coupled BFSDE \eqref{e14} and \eqref{e15}, we are going to decouple it and trying to derive the Nash certainty equivalence (NCE) system satisfied by the decentralized control policy. Then we have
\begin{lemma}\label{l3}
Suppose $P(\cdot)$ is the unique solution of the following Riccati equation system
\begin{equation}\label{e16}
\left \{
\begin{aligned}
&\dot{P}(t)+2AP(t)-B^2R^{-1}P^2(t)+Q=0,\\
&P(T)=H,
\end{aligned}
\right.
\end{equation}
we obtain a revisionary Hamilton system:
\begin{equation}\label{e17}
    \left\{
    \begin{aligned}
    &d\hat{x}_0(t)= \big[A_0\hat{x}_0(t)-B_0^2R_0^{-1}p_0(t)+C_0\hat{z}_0(t)\big]dt+\hat{z}_0(t)dW_{0}(t),\\
    &d\bar{x}(t)=\big[\big(A+D-B^2R^{-1}P(t)\big)\bar{x}(t)-B^2R^{-1}k(t)+\alpha\hat{x}_0(t)\big]dt,\\
    &dk(t)=\big[\big(-A+B^2R^{-1}P(t)\big)k(t)+\big(Q-DP(t)\big)\bar{x}(t)-\alpha P(t)\hat{x}_0(t)\big]dt+\beta_0(t)dW_0(t),\\
    &dp_0(t)=\big[-A_0p_0(t)-Q_0(\hat{x}_0(t)-\bar{x}(t))-\alpha p(t)+\alpha P(t)q(t)\big]dt-C_0p_0(t)dW_0(t),\\
    &dp(t)=\big[-\big(A+D-B^2R^{-1}P(t)\big)p(t)+Q_0(\hat{x}_0(t)-\bar{x}(t))-\big(Q-DP(t)\big)q(t)\big]dt\\
    &\qquad\qquad+\bar{\beta}(t)dW_0(t),\\
    &dq(t)=\big[\big(A-B^2R^{-1}P(t)\big)q(t)+B^2R^{-1}p(t)\big]dt,\\
    &\hat{x}_{0}(T)= \xi,\  \bar{x}(0)=x,\  k(T)=0,\   p_0(0)=-H_0\hat{x}_0(0),\  p(T)=0, \ q(0)=0
    \end{aligned}
    \right.
\end{equation}
which is a triple FBSDE (TFBSDE for short).
\end{lemma}
{\it Proof.} Suppose $$p_i(t)=P_i(t)\hat{x}_i(t)+f_i(t),\ \ 1\leq i\leq N$$ where $P_i(\cdot),f_i(\cdot)$ are to be determined. The terminal condition $p_i(T)= H\hat{x}_i(T)$ implies that $$P_i(T)=H,f_i(T)=0.$$
Applying It\^o's formula to $P_i(t)\hat{x}_i(t)+f_i(t)$, we have
\begin{equation}\nonumber
\begin{aligned}
dp_i(t)=&\big[\dot{P}_i(t)+AP_i(t)-B^2R^{-1}P^2_i(t)\big]\hat{x}_i(t)dt\\
&+\big[DP_i(t)\bar{x}(t)-B^2R^{-1}P_i(t)f_i(t)+\alpha P_i(t)\hat{x}_0(t)\big]dt+df_i(t)+\sigma P_i(t)dW_i(t).
\end{aligned}
\end{equation}
Comparing the coefficients with \eqref{e14}, we get $\beta_i(t)=\sigma P_i(t),$
\begin{equation}\label{e18}
\left \{
\begin{aligned}
&\dot{P}_i(t)+2AP_i(t)-B^2R^{-1}P^2_i(t)+Q=0,\\
&P_i(T)=H
\end{aligned}
\right.
\end{equation}
and
\begin{equation}\label{e19}
\left \{
\begin{aligned}
df_i(t)=&\big[\big(-A+B^2R^{-1}P_i(t)\big)f_i(t)+\big(Q-DP_i(t)\big)\bar{x}(t)-\alpha P_i(t)\hat{x}_0(t)\big]dt\\
&+\beta_0(t)dW_0(t),\\
f_i(T)=&0.
\end{aligned}
\right.
\end{equation}
Noting that Riccati equation \eqref{e18} is symmetric, it is known that \eqref{e18} admits a unique nonnegative bounded solution $P_i(\cdot)$ (see [ma,yong]). Further we get that $P_1(\cdot)=P_2(\cdot)=\cdots=P_N(\cdot):=P(\cdot)$. Then \eqref{e18} coincides with \eqref{e16}. Besides, if the linear BSDE \eqref{e19} admits a unique solution $f_i(\cdot)\in L^{2}_{\mathcal{F}^{w_0}_t}(0, T; \mathbb{R})$, we denote $f_i(\cdot):=f(\cdot),i=1,2,\cdots,N.$

Thus, the decentralized feedback strategy for $\mathcal{A}_i,1\leq i\leq N$ is written as
\begin{equation}\label{e20}
u_{i}(t)=-BR^{-1}\big(P(t)x_i(t)+f(t)\big)
\end{equation}
where $x_i(\cdot)$ and the associated $x_0(\cdot)$ are given by \eqref{e2} and \eqref{e1} respectively. Plugging \eqref{e20} into \eqref{e2} implies the centralized closed-loop state:
\begin{equation}\label{e21}
    \left\{
    \begin{aligned}
    dx_{i}(t)=& \Big[\big(A-B^2R^{-1}P(t)\big)x_i(t)-B^2R^{-1}f(t)+Dx^{(N)}(t)+\alpha x_0(t)\Big]dt+\sigma dW_{i}(t),\\
    x_{i}(0)=& x_{i0}.
    \end{aligned}
    \right.
\end{equation}
Taking summation, dividing by $N$ and letting $N\rightarrow+\infty$, we get
\begin{equation}\label{e22}
    \left\{
    \begin{aligned}
    d\bar{x}(t)=& \big[\big(A+D-B^2R^{-1}P(t)\big)\bar{x}(t)-B^2R^{-1}f(t)+\alpha x_0(t)\big]dt,\\
    \bar{x}(0)=& x.
    \end{aligned}
    \right.
\end{equation}
Comparing the coefficients with \eqref{e5}, we have
\begin{equation}\nonumber
    \begin{aligned}
    \bar{A}=A+D-B^2R^{-1}P(\cdot),\ \ \bar{B}=\alpha,\ \ \bar{C}=-B^2R^{-1},\ \ k(\cdot)=f(\cdot).
    \end{aligned}
\end{equation}
Then we obtain
\begin{equation}\nonumber
\left \{
\begin{aligned}
dk(t)=&\big[\big(-A+B^2R^{-1}P(t)\big)k(t)+\big(Q-DP(t)\big)\bar{x}(t)-\alpha P(t)x_0(t)\big]dt+\beta_0(t)dW_0(t),\\
k(T)=&0.
\end{aligned}
\right.
\end{equation}
From \eqref{e5}, it follows that
\begin{equation}\nonumber
    \begin{aligned}
    \tilde{A}=-A+B^2R^{-1}P(\cdot),\ \ \tilde{B}=Q-DP(\cdot),\ \ \tilde{C}=-\alpha P(\cdot).
    \end{aligned}
\end{equation}
Then \eqref{e17} is obtained, which completes the proof.      \hfill  $\Box$

About the wellposedness of \eqref{e17}, we have the following result.
\begin{theorem}\label{t1}
Suppose $B_0\neq0$, then TFBSDE \eqref{e17} is uniquely solvable.
\end{theorem}
%{\it Proof.} For $\Theta=(\bar{x},p_0,q,p,\hat{x}_0,k,\bar{\beta},\hat{z}_0,\beta_0)$, we denote $\mathbb{A}(t,\Theta):=\Big(-\big(A+D-B^2R^{-1}P(t)\big)p+Q_0(\hat{x}_0-\bar{x})-\big(Q-DP(t)\big)q,A_0\hat{x}_0-B_0^2R_0^{-1}p_0+C_0\hat{z}_0,\big(-A+B^2R^{-1}P(t)\big)k+\big(Q-DP(t)\big)\bar{x}-\alpha P(t)\hat{x}_0,\big(A+D-B^2R^{-1}P(t)\big)\bar{x}-B^2R^{-1}k+\alpha\hat{x}_0,-A_0p_0-Q_0(\hat{x}_0-\bar{x})-\alpha p+\alpha P(t)q,\big(A-B^2R^{-1}P(t)\big)q+B^2R^{-1}p,0,-C_0p_0,0\Big).$
%
%Then for any $\Theta^i=(\bar{x}^i,p_0^i,q^i,p^i,\hat{x}_0^i,k^i,\bar{\beta}^i,\hat{z}_0^i,\beta^i_0),i=1,2,$ we have
%\begin{equation}\nonumber
%    \begin{aligned}
%    \langle\mathbb{A}(t,\Theta^1)-\mathbb{A}(t,\Theta^2),\Theta^1-\Theta^2\rangle
%    =&-Q_0[(\bar{x}^1-\bar{x}^2)-(\hat{x}_0^1-\hat{x}_0^2)]^2-B_0^2R_0^{-1}(p_0^1-p_0^2)^2\\
%    %\leq &-B_0^2R_0^{-1}(p_0^1-p_0^2)^2.
%    \end{aligned}
%\end{equation}
%Then the monotonic condition of FBSDE is satisfied, which implies the wellposedness of \eqref{e17}. The proof is complete.
%\hfill   $\Box$
\begin{remark}
\emph{In what follows \eqref{e17} is called the NCE system. By Theorem \ref{t1} we know that there exists a unique 6-tuple solution $(\bar{x}(\cdot),p_0(\cdot),q(\cdot),p(\cdot),\hat{x}_0(\cdot),k(\cdot))$ which can be obtained outline. Thus it is equivalent with the fixed point principle. To our best knowledge, it is the first time to focus on the wellposedness of TFBSDE in large population problems. It is of great feature and meaningful. }
\end{remark}

%%%%%%%%%%%%%%%%%%%%%%%%%%%%%%%%%%%%%%%%%%%%%% Section 4 %%%%%%%%%%%%%%%%%%%%%%%%%%%%%%%%%%%%%%%%%%%%%%%%%%%%%%%%%%%%%%%%%%%%%%%%%%%%%%%%%%%

\section{$\epsilon$-Nash Equilibrium Analysis}
In above sections, we obtained the optimal control $\bar{u}_i(\cdot), 0\le i\le N$ of Problem (\textbf{II}) through the consistency condition system. Now we turn to verify the $\epsilon$-Nash equilibrium of Problem (\textbf{I}). To start, we first present the definition of $\epsilon$-Nash equilibrium.
\begin{definition}\label{d1}
A set of controls $u_k\in \mathcal{U}_k,\ 0\leq k\leq N,$ for $(N+1)$ agents is called to satisfy an $\epsilon$-Nash equilibrium with respect to the costs $J_k,\ 0\leq k\leq N,$ if there exists $\epsilon\geq0$ such that for any fixed $0\leq i\leq N$, we have
\begin{equation}\label{e24}
J_i(u_i,u_{-i})\leq J_i(u'_i,u_{-i})+\epsilon
\end{equation}
when any alternative control $u'_i\in \mathcal{U}_i$ is applied by $\mathcal{A}_i$.
\end{definition}
If $\epsilon=0,$ then Definition \ref{d1} is reduced to the usual Nash equilibrium. Now, we state the main result of this paper and its proof will be given later.
\begin{theorem}\label{t2}
Under \emph{(H1)}, $(\tilde{u}_0,\tilde{u}_1,\tilde{u}_2,\cdots,\tilde{u}_N)$ satisfies the $\epsilon$-Nash equilibrium of \emph{(\textbf{I})}. Here, $\tilde{u}_0$ is given by
\begin{equation}\label{e25}
\tilde{u}_{0}(t)=-B_0R_0^{-1}p_0(t)
\end{equation}
where $p_0(\cdot)$ is obtained outline by \eqref{e17}; while for $1\le i\le N,$ $\tilde{u}_i$ is
\begin{equation}\label{e26}
\tilde{u}_{i}(t)=-BR^{-1}P(t)\tilde{x}_i(t)-BR^{-1}k(t)
\end{equation}
where $\tilde{x}_{i}(\cdot)$, the decentralized state trajectory for $\mathcal{A}_i$, satisfies \eqref{e21}.
\end{theorem}
The proof of above theorem needs several lemmas which are presented later. Denote by $(\tilde{x}_0(\cdot),\tilde{z}_0(\cdot))$ the centralized state trajectory; $(\hat{x}_0(\cdot),\hat{z}_0(\cdot))$ the decentralized one. Applying $\tilde{u}_0(\cdot)$ to $\mathcal{A}_0$ and using the notations given above, it is easy to get that $(\tilde{x}_0(\cdot),\tilde{z}_0(\cdot))\equiv (\hat{x}_0(\cdot),\hat{z}_0(\cdot))$. Further, $(\bar{x}(\cdot),k(\cdot))_{\tilde{x}_0}=(\bar{x}(\cdot),k(\cdot))_{\hat{x}_0}$. The cost functional in centralized case and decentralized one are given by
\begin{equation}\label{e27}
    \begin{aligned}
    J_{0}(\tilde{u}_{0}(\cdot), \tilde{u}_{-0}(\cdot))=&\frac{1}{2} \mathbb{E}\left\{ \int_{0}^{T}\left[Q_0\Big(\tilde{x}_{0}(t)-\tilde{x}^{(N)}(t)\Big)^{2}+R_0\tilde{u}_{0}^{2}(t)\right]dt+H_{0}\tilde{x}_{0}^{2}(0) \right\}
\end{aligned}
\end{equation}
and \begin{equation}\label{e28}
    \begin{aligned}
    \bar{J}_{0}(\bar{u}_{0}(\cdot))=&\frac{1}{2} \mathbb{E}\left\{ \int_{0}^{T}\left[Q_0\Big(\hat{x}_{0}(t)-\bar{x}(t)_{\hat{x}_0}\Big)^{2}+R_0\bar{u}_{0}^{2}(t)\right]dt+H_{0}\hat{x}_{0}^{2}(0) \right\}
\end{aligned}
\end{equation}
respectively. For $\mathcal{A}_i,1\leq i\leq N$, we have the following close-loop system
\begin{equation}\label{e29}
    \left\{
    \begin{aligned}
    d\tilde{x}_{i}(t)=& \big[(A-B^2R^{-1}P(t))\tilde{x}_i(t)-B^2R^{-1}k(t)_{\tilde{x}_0}+D\tilde{x}^{(N)}(t)+\alpha \tilde{x}_0(t)\big]dt+\sigma dW_{i}(t),\\
    \tilde{x}_{i}(0)=& x_{i0}
    \end{aligned}
    \right.
\end{equation}
with the cost functional
\begin{equation}\label{e30}
    \begin{aligned}
    J_{i}(\tilde{u}_{i}(\cdot), \tilde{u}_{-i}(\cdot))=&\frac{1}{2} \mathbb{E}\left\{ \int_{0}^{T}\left[Q\Big(\tilde{x}_{i}(t)-\tilde{x}^{(N)}(t)\Big)^{2}+R\tilde{u}_{i}^{2}(t)\right]dt+H\tilde{x}_{i}^{2}(T) \right\}
    \end{aligned}
\end{equation}
where $\tilde{x}^{(N)}(t)=\frac{1}{N}\sum\limits^{N}_{i=1}\tilde{x}_{i}(t)$. The auxiliary system (of limiting problem) is given by
\begin{equation}\label{e31}
    \left\{
    \begin{aligned}
    d\hat{x}_{i}(t)=& \big[(A-B^2R^{-1}P(t))\hat{x}_i(t)-B^2R^{-1}k(t)_{\hat{x}_0}+D\bar{x}(t)_{\hat{x}_0}+\alpha \hat{x}_0(t)\big]dt+\sigma dW_{i}(t),\\
    \hat{x}_{i}(0)=& x_{i0}
    \end{aligned}
    \right.
\end{equation}
with the cost functional
\begin{equation}\label{e32}
    \begin{aligned}
    J_{i}(\bar{u}_{i}(\cdot))=&\frac{1}{2} \mathbb{E}\left\{ \int_{0}^{T}\left[Q\Big(\hat{x}_{i}(t)-\bar{x}(t)_{\hat{x}_0}\Big)^{2}+R\bar{u}_{i}^{2}(t)\right]dt+H\hat{x}_{i}^{2}(T) \right\}
    \end{aligned}
\end{equation}
where $(\bar{x}(t)_{\hat{x}_0},k(t)_{\hat{x}_0})$ satisfies \eqref{e17}. We have
\begin{lemma}\label{l4}
\begin{flalign}
\sup_{0\leq t\leq T}\mathbb{E}\Big|\tilde{x}^{(N)}(t)-\bar{x}(t)_{\hat{x}_0}\Big|^2=O\Big(\frac{1}{N}\Big),\label{e33}\\
\Big|J_0(\tilde{u}_0,\tilde{u}_{-0})-\bar{J}_0(\bar{u}_0)\Big|=O\Big(\frac{1}{\sqrt{N}}\Big)\label{e34}.
\end{flalign}
\end{lemma}
{\it Proof.} By \eqref{e29}, we have
\begin{equation*}\left\{\begin{aligned}
d\tilde{x}^{(N)}(t)=&\Big[\big(A+D-B^2R^{-1}P(t)\big)\tilde{x}^{(N)}(t)-B^2R^{-1}k(t)_{\tilde{x}_0}+\alpha \tilde{x}_0(t)\Big]dt+\frac{1}{N}\sum_{i=1}^N\sigma dW_{i}(t),\\
     \tilde{x}^{(N)}(0)=& \frac{1}{N}\sum\limits_{i=1}^Nx_{i0}:=x^{(N)}_0.
\end{aligned}\right.
\end{equation*}
Noting that
\begin{equation}\nonumber\begin{aligned}
\mathbb{E}\Big|x^{(N)}_0-x\Big|^2=\mathbb{E}\left|\int_0^t\frac{1}{N}\sum_{i=1}^N\sigma dW_i(s)\right|^2=O\Big(\frac{1}{N}\Big),
\end{aligned}
\end{equation}
by \eqref{e17} and Gronwall's inequality, we get \eqref{e33}.

It is easily got that $\sup\limits_{0\leq t\leq T}\mathbb{E}\big|\hat{x}_{0}(t)-\bar{x}(t)_{\hat{x}_0}\big|^2<+\infty$. Applying Cauchy-Schwarz inequality, we have
\begin{equation}\label{e35}\begin{aligned}
&\sup_{0\leq t\leq T}\mathbb{E}\Big|\big|\tilde{x}_{0}(t)-\tilde{x}^{(N)}(t)\big|^2-\big|\hat{x}_{0}(t)-\bar{x}(t)_{\hat{x}_0}\big|^2\Big|\\
\leq& \sup_{0\leq t\leq T}\mathbb{E}\big|\tilde{x}_{0}(t)-\tilde{x}^{(N)}(t)-\hat{x}_{0}(t)+\bar{x}(t)_{\hat{x}_0}\big|^2\\
&+2\sup_{0\leq t\leq T}\mathbb{E}\Big[\big|\hat{x}_{0}(t)-\bar{x}(t)_{\hat{x}_0}\big|\big|\tilde{x}_{0}(t)-\tilde{x}^{(N)}(t)-\hat{x}_{0}(t)+\bar{x}(t)_{\hat{x}_0}\big|\Big]\\
\leq&\sup_{0\leq t\leq T}\mathbb{E}\big|\tilde{x}_{0}(t)-\hat{x}_{0}(t)-\big(\tilde{x}^{(N)}(t)-\bar{x}(t)_{\hat{x}_0}\big)\big|^2\\
&+2\Big(\sup_{0\leq t\leq T}\mathbb{E}\big|\hat{x}_{0}(t)-\bar{x}(t)_{\hat{x}_0}\big|^2\Big)^{\frac{1}{2}}\Big(\sup_{0\leq t\leq T}\mathbb{E}\big|\tilde{x}_{0}(t)-\hat{x}_{0}(t)-\big(\tilde{x}^{(N)}(t)-\bar{x}(t)_{\hat{x}_0}\big)\big|^2\Big)^{\frac{1}{2}}\\
=&O\Big(\frac{1}{\sqrt{N}}\Big).
\end{aligned}
\end{equation}
In addition, by \eqref{e10} and \eqref{e25}, we have $\tilde{u}_0(\cdot)=\hat{u}_0(\cdot).$ Thus \eqref{e34} is obtained.   \hfill   $\Box$

For minor players, we have the following estimates:
\begin{lemma}\label{l5}
\begin{flalign}
&\sup_{1\leq i\leq N}\left[\sup_{0\leq t\leq T}\mathbb{E}\Big|\tilde{x}_i(t)-\hat{x}_i(t)\Big|^2\right]=O\Big(\frac{1}{N}\Big),\label{e36}\\
&\sup_{1\leq i\leq N}\left[\sup_{0\leq t\leq T}\mathbb{E}\Big|\tilde{u}_i(t)-\bar{u}_i(t)\Big|^2\right]=O\Big(\frac{1}{N}\Big),\label{e37}\\
&\Big|J_i(\tilde{u}_i,\tilde{u}_{-i})-\bar{J}_i(\bar{u}_i)\Big|=O\Big(\frac{1}{\sqrt{N}}\Big),\ \ 1\leq i\leq N. \label{e38}
\end{flalign}
\end{lemma}
{\it Proof.} For $\forall\ 1\leq i\leq N,$ applying Gronwall's inequality, we get \eqref{e36} from \eqref{e33}. \eqref{e37} follows from \eqref{e36} obviously. Using the same technique as \eqref{e35} and noting $\sup\limits_{0\leq t\leq T}\mathbb{E}\big|\hat{x}_{i}(t)-\bar{x}(t)_{\hat{x}_0}\big|^2<+\infty,\sup\limits_{0\leq t\leq T}\mathbb{E}\big|\bar{u}_{i}(t)\big|^2<+\infty,\sup\limits_{0\leq t\leq T}\mathbb{E}\big|\hat{x}_{i}(t)\big|^2<+\infty$, we obtain \eqref{e38}. \hfill  $\Box$

Until now, we have addressed some estimates of states and costs corresponding to control $\tilde{u}_i$ and $\bar{u}_i$,$0\le i\le N$. Next we will focus on the $\epsilon$-Nash equilibrium for (\textbf{I}). Consider a perturbed control $u_0 \in \mathcal{U}_0$ for $\mathcal{A}_0$ and introduce some notations as
\begin{equation}\label{e39}
    \left\{
    \begin{aligned}
    dl_0(t)=& \big[A_0l_{0}(t)+B_0u_{0}(t)+C_0q_0(t)\big]dt+q_0(t)dW_{0}(t),\\
    x_{0}(T)=& \xi
    \end{aligned}
    \right.
\end{equation}
whereas minor players keep the control $\tilde{u}_i,1\leq i\leq N,$ i.e.,
\begin{equation}\label{e40}
    \left\{
    \begin{aligned}
    dl_{i}(t)=& \big[(A-B^2R^{-1}P(t))l_i(t)-B^2R^{-1}k(t)_{l_0}+Dl^{(N)}(t)+\alpha l_0(t)\big]dt+\sigma dW_{i}(t),\\
    l_{i}(0)=& x_{i0}
    \end{aligned}
    \right.
\end{equation}
where $l^{(N)}(t)=\frac{1}{N}\sum\limits_{k=1}^Nl_k(t)$; $k(t)_{l_0}$ associated with $l_0$ satisfying
\begin{equation}\label{e401}
    \left\{
    \begin{aligned}
    &dk(t)_{l_0}=\big[\big(-A+B^2R^{-1}P(t)\big)k(t)_{l_0}+\big(Q-DP(t)\big)\bar{x}(t)_{l_0}-\alpha P(t)l_0(t)\big]dt\\
    &\qquad\qquad+\beta_0(t)dW_0(t),\\
    &d\bar{x}(t)_{l_0}=\big[\big(A+D-B^2R^{-1}P(t)\big)\bar{x}(t)_{l_0}-B^2R^{-1}k(t)_{l_0}+\alpha l_0(t)\big]dt,\\
    &k(T)_{l_0}=0,\  \bar{x}(0)_{l_0}=x.
    \end{aligned}
    \right.
\end{equation}
And for any fixed $i$, $1\le i\le N$, consider a perturbed control $u_i \in \mathcal{U}_i$ for $\mathcal{A}_i$, whereas major and other minor players keep the control $\tilde{u}_j,0\leq j\leq N,j\neq i.$ Introduce
\begin{equation}\label{e41}
    \left\{
    \begin{aligned}
    dm_{i}(t)=& \big[Am_{i}(t)+Bu_{i}(t)+Dm^{(N)}(t)+\alpha \tilde{x}_0(t)\big]dt+\sigma dW_{i}(t),\\
    m_{i}(0)=& x_{i0}
    \end{aligned}
    \right.
\end{equation} and
\begin{equation}\label{e42}
    \left\{
    \begin{aligned}
    dm_{j}(t)=& \big[(A-B^2R^{-1}P(t))m_j(t)-B^2R^{-1}k(t)_{\tilde{x}_0}+Dm^{(N)}(t)+\alpha \tilde{x}_0(t)\big]dt+\sigma dW_{i}(t),\\
    m_{j}(0)=& x_{j0}
    \end{aligned}
    \right.
\end{equation}
where $m^{(N)}(t)=\frac{1}{N}\sum\limits_{k=1}^Nm_k(t)$; $k(t)_{\tilde{x}_0}$ satisfying \eqref{e17} because of $\tilde{x}_0=\hat{x}_0$.

If $\tilde{u}_j,\ 0\leq j\leq N$ is an $\epsilon$-Nash equilibrium with respect to cost $J_j$, it holds that
$$J_j(\tilde{u}_j,\tilde{u}_{-j})\geq \inf_{u_j\in \mathcal{U}_j}J_i(u_j,\tilde{u}_{-j})\geq J_j(\tilde{u}_j,\tilde{u}_{-j})-\epsilon.$$
Then, when making the perturbation, we just need to consider $u_j\in \mathcal{U}_j$ such that $J_j(u_j,\tilde{u}_{-j})\leq J_j(\tilde{u}_j,\tilde{u}_{-j}),$ which implies
\begin{equation}\nonumber\begin{aligned}
\frac{1}{2}\mathbb{E}\int_0^TRu_j^2(t)dt\leq J_j(u_i,\tilde{u}_{-j})\leq J_j(\tilde{u}_j,\tilde{u}_{-j})=\bar{J}_j(\bar{u}_j)+O\Big(\frac{1}{\sqrt{N}}\Big).
\end{aligned}
\end{equation}
In the limiting cost functional $\bar{J}_j$, by the optimality of $(\bar{x}_j,\bar{u}_j)$, we get that $(\bar{x}_j,\bar{u}_j)$ is $L^2$-bounded. Then we obtain the boundedness of $\bar{J}_j(\bar{u}_j)$, i.e.,
\begin{equation}\label{e43}
    \mathbb{E}\int_0^Tu_j^2(t)dt\leq C,\ \ 0\leq j\leq N
\end{equation}where $C$ is a positive constant and independent of $N$. Then we have the following proposition.
\begin{proposition}\label{pr1}
 $\sup\limits_{0\leq t\leq T}\mathbb{E}\big|l_0(t)\big|^2$, $\sup\limits_{1\le k\le N}\Bigg[\sup\limits_{0\leq t\leq T}\mathbb{E}\big|l_k(t)\big|^2 \Bigg]$, $\sup\limits_{1\le k\le N}\Bigg[\sup\limits_{0\leq t\leq T}\mathbb{E}\big|m_k(t)\big|^2 \Bigg]$ are bounded.
\end{proposition}
{\it Proof.} By \eqref{e43}, using the usual technique of BSDE, we get the boundedness of $\sup\limits_{0\leq t\leq T}\mathbb{E}\big|l_0(t)\big|^2$. It follows from \eqref{e40} that
\begin{equation}\nonumber\begin{aligned}
\mathbb{E}\Big[\sum_{k=1}^N|l_k(t)|^2\Big]\leq & C_1\Bigg\{\mathbb{E}\Big[\sum_{k=1}^N|x_{k0}|^2\Big]+\mathbb{E}\int_0^t\Big[\sum_{k=1}^N|l_k(s)_{l_0}|^2+N|k(s)|^2+N|l_0(s)|^2\Big]ds\\
&+\sum_{k=1}^N\mathbb{E}\Big|\int_0^t\sigma dW_k(s)\Big|^2\Bigg\}.
\end{aligned}
\end{equation}
By \eqref{e41} and \eqref{e42}, it holds that
\begin{equation}\nonumber\begin{aligned}
\mathbb{E}\Big[\sum_{k=1}^N|m_k(t)|^2\Big]\leq &C_2\Bigg\{\mathbb{E}\Big[\sum_{k=1}^N|x_{k0}|^2\Big]+\mathbb{E}\int_0^t\Big[\sum_{k=1}^N|m_k(s)|^2+|u_i(s)|^2+\sum_{k=1,k\neq i}^N|\tilde{u}_k(s)|^2\\
&+N|\tilde{x}_0(s)|^2\Big]ds+\sum_{k=1}^N\mathbb{E}\Big|\int_0^t\sigma dW_k(s)\Big|^2\Bigg\}.
\end{aligned}
\end{equation}
Here, $C_1,C_2$ are both positive constants. Since $\sup\limits_{0\leq t\leq T}\mathbb{E}\big|l_0(t)\big|^2$ is bounded, we get the boundedness of $\sup\limits_{0\leq t\leq T}\mathbb{E}\big|k(t)_{l_0}\big|^2$ by \eqref{e401}. It follows from \eqref{e43} that $\mathbb{E}|u_i(t)|^2$ is bounded. Besides, the optimal controls $\tilde{u}_k(t),k\neq i$ is $L^2$-bounded. Then by Gronwall's inequality, it follows that
\begin{equation}\nonumber
\sup_{0\leq t\leq T}\mathbb{E}\left[\sum_{k=1}^N|l_k(t)|^2\right]=\sup_{0\leq t\leq T}\mathbb{E}\left[\sum_{k=1}^N|m_k(t)|^2\right]=O(N).
\end{equation}
Hence for any $1\leq k\leq N,$ $\sup\limits_{0\leq t\leq T}\mathbb{E}|l_k(t)|^2$ and $\sup\limits_{0\leq t\leq T}\mathbb{E}|m_k(t)|^2$ are bounded.    \hfill  $\Box$

Correspondingly, the system for agent $\mathcal{A}_0$ under control $u_0$ in \textbf{(II)} is as follows
\begin{equation}\label{e44}
    \left\{
    \begin{aligned}
    dl'_0(t)=& \big[A_0l'_{0}(t)+B_0u_{0}(t)+C_0q'_0(t)\big]dt+q'_0(t)dW_{0}(t),\\
    x'_{0}(T)=& \xi
    \end{aligned}
    \right.
\end{equation}
and for agent $\mathcal{A}_i,1\leq i\leq N$,
\begin{equation}\label{e45}
    \left\{
    \begin{aligned}
    d\hat{l}_{i}(t)=& \big[(A-B^2R^{-1}P(t))\hat{l}_i(t)-B^2R^{-1}k(t)_{l'_0}+D\bar{x}(t)_{l'_0}+\alpha l'_0(t)\big]dt+\sigma dW_{i}(t),\\
    \hat{l}_{i}(0)=& x_{i0}.
    \end{aligned}
    \right.
\end{equation}
where $(k(t)_{l'_0},\bar{x}(t)_{l'_0})$ associated with $l'_0$ satisfying
\begin{equation}\label{e46}
    \left\{
    \begin{aligned}
    &dk(t)_{l'_0}=\big[\big(-A+B^2R^{-1}P(t)\big)k(t)_{l'_0}+\big(Q-DP(t)\big)\bar{x}(t)_{l'_0}-\alpha P(t)l'_0(t)\big]dt\\
     &\qquad\qquad+\beta_0(t)dW_0(t),\\
    &d\bar{x}(t)_{l'_0}=\big[\big(A+D-B^2R^{-1}P(t)\big)\bar{x}(t)_{l'_0}-B^2R^{-1}k(t)_{l'_0}+\alpha l'_0(t)\big]dt,\\
    &k(T)_{l'_0}=0,\  \bar{x}(0)_{l'_0}=x.
    \end{aligned}
    \right.
\end{equation}
 Then we have
\begin{lemma}\label{l6}
\begin{flalign}
\sup_{0\leq t\leq T}\mathbb{E}\Big|l^{(N)}(t)-\bar{x}(t)_{l'_0}\Big|^2=O\Big(\frac{1}{N}\Big),\label{e47}\\
\Big|J_0(u_0,\tilde{u}_{-0})-\bar{J}_0(u_0)\Big|=O\Big(\frac{1}{\sqrt{N}}\Big)\label{e48}.
\end{flalign}
\end{lemma}
{\it Proof.} From \eqref{e39} and \eqref{e44}, by the existence and uniqueness of BSDE, we have $(l'_0,q'_0)=(l_0,q_0)$. Further, noting FBSDE \eqref{e401} and \eqref{e46}, we get $(k(t)_{l'_0},\bar{x}(t)_{l'_0})=(k(t)_{l_0},\bar{x}(t)_{l_0})$.

From \eqref{e40}, it follows that
\begin{equation*}\left\{\begin{aligned}
dl^{(N)}(t)=&\Big[\big(A+D-B^2R^{-1}P(t)\big)l^{(N)}(t)-B^2R^{-1}k(t)_{l_0}+\alpha l_0(t)\Big]dt+\frac{1}{N}\sum_{i=1}^N\sigma dW_{i}(t),\\
     l^{(N)}(0)=& x^{(N)}_0.
\end{aligned}\right.
\end{equation*}
Noting \eqref{e46} and
\begin{equation}\nonumber\begin{aligned}
\mathbb{E}\Big|x^{(N)}_0-x_0\Big|^2=\mathbb{E}\left|\int_0^t\frac{1}{N}\sum_{i=1}^N\sigma dW_i(s)\right|^2=O\Big(\frac{1}{N}\Big),
\end{aligned}
\end{equation}
and applying Gronwall's inequality, we get \eqref{e47}. Using the same technique as Lemma \ref{l4} and noting $\sup\limits_{0\leq t\leq T}\mathbb{E}\big|l'_{0}(t)-\bar{x}(t)_{l'_0}\big|^2<+\infty$, we obtain \eqref{e48}.    \hfill  $\Box$

Now, we will focus on the difference for the perturbed control and optimal control of minor players. Given the system for agent $\mathcal{A}_i$ under control $u_i$ in \textbf{(II)} is
\begin{equation}\label{e49}
    \left\{
    \begin{aligned}
    dm'_{i}(t)=& \big[Am'_{i}(t)+Bu_{i}(t)+D\bar{x}(t)_{\hat{x}_0}+\alpha \hat{x}_0(t)\big]dt+\sigma dW_{i}(t),\\
    m'_{i}(0)=& x_{i0}
    \end{aligned}
    \right.
\end{equation} and for agent $\mathcal{A}_j,j\neq i$,
\begin{equation}\label{e50}
    \left\{
    \begin{aligned}
    d\hat{m}_{j}(t)=& \big[(A-B^2R^{-1}P(t))\hat{m}_j(t)-B^2R^{-1}k(t)_{\hat{x}_0}+D\bar{x}(t)_{\hat{x}_0}+\alpha \hat{x}_0(t)\big]dt+\sigma dW_{i}(t),\\
    \hat{m}_{j}(0)=& x_{j0}
    \end{aligned}
    \right.
\end{equation}
where $(\bar{x}(t)_{\tilde{x}_0},k(t)_{\tilde{x}_0})$ satisfying \eqref{e17} because of $\tilde{x}_0=\hat{x}_0$.

In order to give necessary estimates in (\textbf{I}) and (\textbf{II}), we also introduce some intermediate states as
\begin{equation}\label{e51}\left\{\begin{aligned}
d\check{m}_i(t)=&\left[A\check{m}_i(t)+Bu_i(t)+\frac{N-1}{N}D \check{m}^{(N-1)}(t)+\alpha \tilde{x}_0(t)\right]dt +\sigma dW_i(t),\\
\check{m}_i(0)=&x_{i0}
\end{aligned}\right.
\end{equation}and for $j\neq i$,
\begin{equation}\label{e52}\left\{\begin{aligned}
d\check{m}_j(t)=&\left[\Big(A-B^2R^{-1}P(t)\Big)\check{m}_{j}(t)-B^2R^{-1}k(t)_{\tilde{x}_0}+\frac{N-1}{N}D \check{m}^{(N-1)}(t)+\alpha \tilde{x}_0(t)\right]dt\\
&+\sigma dW_{j}(t),\\
\check{m}_j(0)=&x_{j0}
\end{aligned}\right.
\end{equation}
where $\check{m}^{(N-1)}(t)=\frac{1}{N-1}\sum\limits_{j=1,j\neq i}^N\check{m}_j(t)$.

Denoting $m^{(N-1)}(t)=\frac{1}{N-1}\sum\limits_{j=1,j\neq i}^Nm_j(t)$, by \eqref{e42} and \eqref{e52}, we get
\begin{equation}\label{e53}\left\{\begin{aligned}
dm^{(N-1)}(t)=&\left[\big(A-B^2R^{-1}P(t)+\frac{N-1}{N}D\big)m^{(N-1)}(t)-B^2R^{-1}k(t)_{\tilde{x}_0}+\alpha \tilde{x}_0(t)+\frac{D}{N}m_i(t)\right]dt\\
&+\frac{1}{N-1}\sum_{j=1,j\neq i}^N\sigma dW_j(t),\\
m^{(N-1)}(0)=&\frac{1}{N-1}\sum\limits_{j=1,j\neq i}^Nx_{j0}:=x^{(N-1)}_0
\end{aligned}\right.
\end{equation}and
\begin{equation}\label{e54}\left\{\begin{aligned}
d\check{m}^{(N-1)}(t)=&\left[\big(A-B^2R^{-1}P(t)+\frac{N-1}{N}D\big)\check{m}^{(N-1)}(t)-B^2R^{-1}k(t)_{\tilde{x}_0}+\alpha \tilde{x}_0(t)\right]dt\\
&+\frac{1}{N-1}\sum_{j=1,j\neq i}^N\sigma dW_j(t),\\
\check{m}^{(N-1)}(0)=&x^{(N-1)}_0.
\end{aligned}\right.
\end{equation}
We have the following estimates on these states.
\begin{proposition}\label{Pr2}
\begin{align}
\label{e55}
&\sup_{0\leq t\leq T}\mathbb{E}\Big|m^{(N-1)}(t)-\check{m}^{(N-1)}(t)\Big|^2=O\Big(\frac{1}{N^2}\Big),\\
\label{e56}
&\sup_{0\leq t\leq T}\mathbb{E}\Big|m^{(N)}(t)-m^{(N-1)}(t)\Big|^2=O\Big(\frac{1}{N}\Big),\\
\label{e57}
&\sup_{0\leq t\leq T}\mathbb{E}\Big|\check{m}^{(N-1)}(t)-\bar{x}(t)_{\hat{x}_0}\Big|^2=O\Big(\frac{1}{N}\Big).
\end{align}
\end{proposition}
{\it Proof.} From \eqref{e53}-\eqref{e54}, applying Proposition \ref{pr1} and Gronwall's inequality, the assertion \eqref{e55} holds. \eqref{e56} follows from assumption (H1) and the $L^2$-boundness of controls $u_i(\cdot)$ and $\tilde{u}_j(\cdot),j\neq i.$ From \eqref{e54} and \eqref{e17}, noting $(\bar{x}(t)_{\tilde{x}_0},k(t)_{\tilde{x}_0},\tilde{x}_0)=(\bar{x}(t)_{\hat{x}_0},k(t)_{\hat{x}_0},\hat{x}_0)$, we get
\begin{equation}\nonumber\left\{\begin{aligned}
d\Big(\check{m}^{(N-1)}(t)-\bar{x}(t)_{\hat{x}_0}\Big)=&\Big[\frac{N-1}{N}D\big(\check{m}^{(N-1)}(t)-\bar{x}(t)_{\hat{x}_0}\big)-\frac{D}{N}\bar{x}(t)_{\hat{x}_0}\Big]dt +\frac{1}{N-1}\sum_{j=1,j\neq i}^N\sigma dW_j(t),\\
\check{m}^{(N-1)}(0)-\bar{x}(0)_{\hat{x}_0}=&x^{(N-1)}_0-x.
\end{aligned}\right.
\end{equation}
Thus, \eqref{e57} is obtained.     \hfill  $\Box$

Based on Proposition \ref{Pr2}, we obtain more direct estimates to prove Theorem \ref{t2}.
\begin{lemma}\label{l7} For fixed $i,1\leq i\leq N$, we have
\begin{align}\label{e58}
&\sup_{0\leq t\leq T}\mathbb{E}\Big|m^{(N)}(t)-\bar{x}(t)_{\hat{x}_0}\Big|^2=O\Big(\frac{1}{N}\Big),\\ \label{e59}
&\sup_{0\leq t\leq T}\mathbb{E}\Big|m_i(t)-m'_i(t)\Big|^2=O\Big(\frac{1}{N}\Big),\\ \label{e60}
&\Big|J_i(u_i,\tilde{u}_{-i})-\bar{J}_i(u_i)\Big|=O\Big(\frac{1}{\sqrt{N}}\Big).
\end{align}
\end{lemma}
{\it Proof.} \eqref{e58} follows from Proposition \ref{Pr2} directly. From \eqref{e41} and \eqref{e49}, we get \eqref{e59} by applying the result of \eqref{e58}.
Further, we have
\begin{equation}\nonumber\begin{aligned}
\sup_{0\leq t\leq T}\mathbb{E}\Big||m_i(t)|^2-|m'_i(t)|^2\Big|=O\Big(\frac{1}{\sqrt{N}}\Big).
\end{aligned}
\end{equation}
In addition,
\begin{equation}\nonumber\begin{aligned}
&\sup_{0\leq t\leq T}\mathbb{E}\left|\big(m_{i}(t)-m^{(N)}(t)\big)^2-\big(m'_{i}(t)-\bar{x}(t)_{\hat{x}_0}\big)^{2}\right|\\
\leq&\sup_{0\leq t\leq T}\mathbb{E}\Big|m_{i}(t)-m'_{i}(t)-\big(m^{(N)}(t)-\bar{x}(t)_{\hat{x}_0}\big)\Big|^2\\
&+2\Big(\sup_{0\leq t\leq T}\mathbb{E}\big|m'_{i}(t)-\bar{x}(t)_{\hat{x}_0}\big|^2\Big)^{\frac{1}{2}}\Big(\sup_{0\leq t\leq T}\mathbb{E}\big|m_{i}(t)-m'_{i}(t)-\big(m^{(N)}(t)-\bar{x}(t)_{\hat{x}_0}\big)\big|^2\Big)^{\frac{1}{2}}\\
=&O\Big(\frac{1}{\sqrt{N}}\Big).
\end{aligned}
\end{equation}
%\Big|J_i(\tilde{u}_i,\tilde{u}_{-i})-\bar{J}_i(\bar{u}_i)\Big|\\
%\leq& \quad\frac{1}{2}\mathbb{E}\int_0^T \Big[ Q\Big|\Big(\tilde{x}_{i}(t)-\big(S\tilde{x}^{(N)}(t)+\eta\big)\Big)^2-\Big(\hat{x}_{i}(t)-\big(S\bar{x}(t)+\eta\big)\Big)^{2}\Big|+R\Big|\tilde{u}_i^2(t)-\bar{u}_i^2(t)\Big|\Big]dt\\
%&+ \frac{1}{2}N_0\mathbb{E}\Big|\tilde{y}_{i}^{2}(0)-\hat{y}_{i}^{2}(0)\Big|\\
%=& O\Big(\frac{1}{\sqrt{N}}\Big),
Then we have \begin{equation}\nonumber\begin{aligned}
&\Big|J_i(u_i,\tilde{u}_{-i})-\bar{J}_i(u_i)\Big|\\
\leq&\quad\frac{1}{2}\mathbb{E}\int_0^T Q\Big|\big(m_{i}(t)-m^{(N)}(t)\big)^2-\big(m'_{i}(t)-\bar{x}(t)_{\hat{x}_0}\big)^{2}\Big|dt\\
&+ \frac{1}{2}H\mathbb{E}\Big|m_{i}^{2}(T)-\big(m'_{i}(T)\big)^{2}\Big|\\
=& O\Big(\frac{1}{\sqrt{N}}\Big)
\end{aligned}
\end{equation}
which implies \eqref{e60}.           \hfill  $\Box$

\emph{Proof of Theorem} \ref{t2}: Now, we consider the $\epsilon$-Nash equilibrium for $\mathcal{A}_0$ and $\mathcal{A}_i,1\leq i\leq N$. Combining \eqref{e34} and \eqref{e48}, we have
\begin{equation}\nonumber\begin{aligned}
J_0(\tilde{u}_0,\tilde{u}_{-0})&=\bar{J}_0(\bar{u}_0)+O\Big(\frac{1}{\sqrt{N}}\Big)\\
&\leq \bar{J}_0(u_0)+O\Big(\frac{1}{\sqrt{N}}\Big)\\
&=J_0(u_0,\tilde{u}_{-0})+O\Big(\frac{1}{\sqrt{N}}\Big).
\end{aligned}
\end{equation}
While from \eqref{e38} and \eqref{e60}, we get
\begin{equation}\nonumber\begin{aligned}
J_i(\tilde{u}_i,\tilde{u}_{-i})&=\bar{J}_i(\bar{u}_i)+O\Big(\frac{1}{\sqrt{N}}\Big)\\
&\leq \bar{J}_i(u_i)+O\Big(\frac{1}{\sqrt{N}}\Big)\\
&=J_i(u_i,\tilde{u}_{-i})+O\Big(\frac{1}{\sqrt{N}}\Big).
\end{aligned}
\end{equation}
Thus, Theorem \ref{t2} follows by taking $\epsilon=O\Big(\frac{1}{\sqrt{N}}\Big)$.

\section{Conclusion}

In this paper, the dynamics of major player is given by some BSDE, while dynamics of minor players are described by SDEs. The BFSDE system is established in which a large number of negligible agents are coupled in their dynamics via state average. Some auxiliary mean-field SDEs (MFSDEs) and a $3\times2$ mixed FBSDE system are considered and analyzed instead of involving the fixed-point analysis as in other mean-field games. The decentralized strategies are derived, which are also shown to satisfy the $\epsilon$-Nash equilibrium property. In the future, one possible direction is that state-average appears in dynamics of major player, which may bring lots of trouble in proving the $\epsilon$-Nash equilibrium. Another direction is that dynamics of minor players are backward and the consistent condition analysis may be more complicated. Numerical computation and other applications in finance will also be investigated in the future work.

%\end{multicols}%%%%%%%%%%%%%%%%%%%%%%%%%%%%%%%%%%%%%%%%%%%%%%%%%%%%%%%%%%%%%%%%%%%%%%%%%%%%%%%%%%%%%%%%%%%%%%%%%%%%%%%%%%%%%%%%%%%%%%

\end{document}